\documentclass{article}

\usepackage {amsfonts}
\usepackage {amssymb}
\usepackage {amsmath}
\usepackage{amsbsy}
\usepackage {hyperref}
\usepackage{bbm} 
\usepackage{graphicx}
\usepackage {color}
\usepackage{amsthm}

\newtheorem{Ps}{Proposition}

\newtheorem{T}{Theorem}
\newtheorem{Lm}{Lemma}
\newtheorem{D}{Definition}

\author{Rapha\"el Lachi\`eze-Rey\thanks{Univ. Lille 1, raphael.lachieze-rey@math.univ-lille1.fr}}
\title{Strong mixing property for STIT tessellations}

\newcommand{\Sd}{\mathcal{S}^{d-1}}
\newcommand{\A}{\forall}
\newcommand{\Ex}{\exists}
\def \d{\textrm{d}}

\begin{document}

\maketitle

\paragraph{abstract}
The so-called STIT tessellations form the class of homogeneous (spatially stationary) tessellations of $\mathbb{R}^d$ which are stable under the nesting/iteration operation. In this paper, we establish the strong mixing property for these tessellations and give the optimal form of the rate of decay for the quantity $|\mathbb{P}({A}\cap Y=\emptyset,T_h B \cap Y=\emptyset)-\mathbb{P}({A}\cap Y=\emptyset)\mathbb{P}({B}\cap Y=\emptyset)|$ when $A$ and $B$ are two compact sets, $h$ a vector of $\mathbb{R}^d$, $T_{h}$ the corresponding translation operator and $Y$ a STIT Tessellation.

~\\

\paragraph{keywords}Stochastic geometry; random tessellations; iteration of tessellations; stability; space ergodicity; strong mixing property; rate of mixing
~\\

\section{Introduction and notations}

Random tessellations, or mosaics, form an important class of objects of stochastic geometry. They have proven to be a useful tool for modeling geometrical structures appearing in biology, geology, medical sciences, and many others. Poisson hyperplane tessellation, Poisson Voronoi tessellation, and its dual tessellation Poisson Delaunay tessellation, are the most celebrated and tractable models investigated until now, and all are defined using a Poisson point process on an appropriate space. In the 1980's, Ambartzumian had the idea of an operation of iteration among the mosaics of $\mathbb{R}^d$, namely the iteration, also called nesting operation. When regarded in $\mathbb{R}^2$, the nesting operation leaves stable the class of T-noded tessellations, while the superposition leaves stable the class of X-noded tessellations, such as the Poisson line tessellation. Nagel, Mecke and Weiss (\cite{MNW1},\cite{MNW2},\cite{NW1}) have introduced the model of the STIT tessellation, following the idea of Cowan \cite{C84}, which is a very natural object regarding this operation, because it possesses the feature of being STable under ITeration. This object can be used as a model for crack patterns, such as those seen on old pottery or during earthquakes.\\
The tessellation is homogeneous, i.e space stationary, and a proper choice of paramaters can make it isotropic, which yields a very interesting model. Many geometrical features of STIT tessellations have been investigated, including moments of variables related to typical faces of the tessellation. Cowan, in \cite{C78} and \cite{C80}, emphasized the importance of ergodic properties for random tessellations. In this paper we establish that all STIT tessellations possess the strong mixing property, which implies its ergodicity. Namely, if $A$ and $B$ are two Borel sets, and $Y$ is the closed set of boundaries of the cells of a STIT tesselation, then

\begin{eqnarray}
\label{EquMix}
\mathbb{P}(Y \cap A = \emptyset,Y \cap T_h B = \emptyset)-\mathbb{P}(Y \cap A = \emptyset)\mathbb{P}(Y \cap B = \emptyset) \xrightarrow[||h|| \to \infty]{} 0,
\end{eqnarray}
where $T_h B$ is the set $B$ translated according to the vector $h$.
 Moreover, when $A$ and $B$ are compact sets, we are able to give optimal rate of decay, of order $\frac{1}{||h||}$.\\
After introducing the notations, we will give a brief decription of the construction and the properties of STIT tessellations, and we will finally establish its mixing properties. The last section is devoted to the proofs.\\

We introduce the following notations:
\begin{itemize}
\item  $E=\mathbb{R}^d$ is the $d$-dimensional euclidean space.
\item  For any set $\Omega$ equipped with an unambiguous topology, $\mathcal{B}(\Omega)$ is the $\sigma$-algebra of Borel sets of $\Omega$.
\item $\mathcal{F}$ is the class of closed sets of $E$.
\item $\mathcal{G}$ is the class of open sets of $E$.
\item $\mathcal{K}$ is the class of compact sets of $E$.
\item $\mathcal{H}$ is the set of affine hyperplanes of $E$.
\item $\Sd$ is the unit sphere of $E$.
\item $\mathbb{R}^+$ is the set of all non-negative numbers.
\end{itemize}

For $A$ a subset of $E$:
\begin{itemize}
\item $\mathop{\mathrm{int}}(A)$ denotes the interior of $A$.
\item $\partial(A)$ denotes the boundary of $A$.
\item $\mathop{\mathrm{conv}}(A)$ denotes the convex hull of $A$.
\item $\mathop{\mathrm{span}}(A)$ is the smallest linear subspace of $E$ containing $A$.
\end{itemize}
For a scalar $\lambda \geq 0$ and a vector $h \in \mathbb{R}^d$, the scaling and translation operators (resp. $x \to \lambda x$ and $T_{h}:x \to x+h$) are naturally lifted to $\mathcal{F}$, with the same notations.\\
Since we are interested in hitting and missing probabilities here, we introduce the corresponding families.
\begin{D}
For $A \in \mathcal{B}(E)$,
\begin{eqnarray*}
\mathcal{F}_{A}=\{C \in \mathcal{F}~|~C \cap A \neq \emptyset \},\\
\mathcal{F}^{A}=\{C \in \mathcal{F}~|~C \cap A = \emptyset \}.
\end{eqnarray*}
\end{D}
\begin{D}
The $\textbf{toplology of closed convergence}$ on $\mathcal{F}$, or $\textit{Fell}$ topology, is the topology generated by the $\mathcal{F}^K,K \in \mathcal{K}$ and the  $\mathcal{F}_O,O \in \mathcal{G}$.
\end{D}
For any set $\mathcal{P} \subset \mathcal{F}$, we can define the Fell topology by restriction and the corresponding Borel $\sigma$-algebra $\mathcal{B}(\mathcal{P})$.
 We know that the $\sigma$-algebra $\mathcal{B}(\mathcal{F})$  is generated by the $\mathcal{F}^K,K \in \mathcal{K}$ (see \cite{M}).
Let us give definitions for tessellations and random tessellations.

\begin{D}
We define a tessellation of  $\mathbb{R}^d$ as a countable set $R$ of convex polytopes of $\mathbb{R}^d$ which satisfy
\begin{description}
\item[(i)] $\A ~C \in R, \mathop{\textrm{int}}(C) \neq \emptyset$.
\item[(ii)] $ \bigcup_{C \in R}C=\mathbb{R}^d.$
\item[(iii)] $ \A ~C,C' \in R, ~C \neq C' \Rightarrow \mathop{\textrm{int}}(C) \bigcap \mathop{\textrm{int}}(C') = \emptyset$.
\item[(iv)] $\A ~K \in \mathcal{K}, ~\mathop{\textrm{card}}(\mathcal{F}_{K} \bigcap R) < \infty$.
\end{description}
$\mathcal{M}_{C}=\{R \in \mathcal{P}(\mathcal{K})~|~R ~\textit{satisfies}~ (i),(ii),(iii),(iv)\}$.
\end{D}~\\

We define the application $\mathcal{Y}$, from $\mathcal{M}_c$ to $\mathcal{F}$, by 

\begin{eqnarray*}
\mathcal{Y}(R)= \bigcup_{C \in R} \partial C,
\end{eqnarray*}
 the closed set formed by the boundaries of all cells of the tessellation.

We can define the reciprocal application: For $Y$ a  closed set, we define $\mathcal{R}(Y)$ as the set of the topological closures of all connected components of $Y's$ complementary.\\

\begin{D}
Let $\mathcal{M}$ be the class of tesselations as closed sets of $\mathbb{R}^d$:
\begin{eqnarray*}
\mathcal{M}=\{Y \in \mathcal{F}~|~\mathcal{R}(Y) \in \mathcal{M}_{C}\}.
\end{eqnarray*}~\\

\end{D}
 Let $\mathcal{B}(\mathcal{M}_c)$ be the $\sigma$-algebra on $\mathcal{M}_c$ that makes $$\mathcal{Y}:\left(\mathcal{M},\mathcal{B}(\mathcal{M})\right)\to (\mathcal{M}_c,\mathcal{B}(\mathcal{M}_c))$$ measurable (where $\mathcal{M}$ is endowed with the restriction of the Fell topology).\\
 Since the applications $\mathcal{Y}$ and $\mathcal{R}$ are measurable reciprocal bijections, we can define a random tessellation as a random variable taking values in either $\mathcal{M}_{C}$ or $\mathcal{M}$, and henceforth we will often jump from one description to another, without bothering with the measurability questions.\\
According to Choquet's theorem (see \cite{M}, \cite{SW}), the law of a random closed set is characterized by its capacity functional, which is defined by:

\begin{D}
Let $Y$ be a random closed set of law $P$ on $\mathcal{B}(\mathcal{F})$. We define\\
\begin{eqnarray*}
\A ~K \in \mathcal{K}, ~  T_Y(K)=P(\mathcal{F}_K).\\
\end{eqnarray*}~\\
\end{D}

If one is given two random closed sets $Y$ and $Z$, one will be able to say that they have the same law if there is equality of their capacity functionals.
\begin{eqnarray*}
Y \stackrel {(d)} {=} Z \Leftrightarrow \A ~K \in \mathcal{K},~ T_Y(K)=T_Z(K).
\end{eqnarray*}

\section{STIT Tessellations}
\label{SecSTIT}
We will first give the decription of the local construction of a STIT tessellation, and then use spatial consistency to extend the definition to the whole space. We will then define the nesting iteration operation and give its connections with STIT tessellations, and finally give the main properties of such tessellations.\\
A STIT tessellation, like the Poisson hyperplane tessellation, is defined from an hyperplane Poisson process. All details oncerning the construction of STIT tessellations can be found in \cite{MNW1}, \cite{MNW2} and \cite{NW1}.

Let us begin with some definitions.
\begin{D}
For $A \in \mathcal{B}$, we set
\begin{eqnarray*}
[A]=\mathcal{H} \cap \mathcal{F}_A
\end{eqnarray*}
the hyperplanes hitting $A$.\\
For $B$ another borelian set, we also define $[A|B]$ as  the set of all hyperplanes $\gamma$ that separate strictly $A$ and $B$, i.e such that $A$ and $B$ are contained in different open half-spaces of $E \setminus \gamma$. \\
We write in short, for $x,y \in \mathbb{R}^d$,
\begin{eqnarray*}
[x]:=[\{x\} ],\\
and~[x|y]:=[\{ x \} |\{y\} ].\\
\end{eqnarray*}
We equip $\mathcal{H}$ with the Borel $\sigma$-algebra induced by Fell topology and let $\Lambda$ be a translation  invariant locally finite measure on $\mathcal{H}$.\\
Given $\gamma \in \mathcal{H} \setminus [0]$, we denote $\gamma^+$ the closed half-space delimited by $\gamma$ not containing $0$, and $\gamma^-$ the other one. Note that $\Lambda([0])=0$ since $\Lambda$ is stationary and locally finite.\\
\end{D}
Let $a$ be a positive number.\\
We identify $\mathcal{H}$ with $\mathbb{R}^+ \times \mathcal{S}^{d-1}$ where an element $(r,u)$ of $\mathbb{R}^+ \times \mathcal{S}^{d-1}$ stands for the hyperplane $\gamma$ at distance $r$ of the origin and normal exterior vector $u$ (i.e normal vector directed towards $\gamma^+$).\\
Since $\Lambda$ is assumed to be stationary, we can write, along with this identification
\begin{eqnarray*}
\Lambda = \lambda^+ \otimes \nu,
\end{eqnarray*}
where $\lambda^+$ is the Lebesgue measure restricted to $\mathbb{R}^+$ and $\nu$ is a finite measure on $\mathcal{S}^{d-1}$.\\
We assume that $\nu$ has the following property:
\begin{eqnarray}
\label{ConNu2}
\textrm{span}(\textrm{supp}(\nu))=E.
\end{eqnarray}
The latter property ensures that  tessellations generated by the hyperplane processes with distribution $\Lambda$ will a.s have compact cells.\\
On a compact  window $W \in \mathcal{K}$ with non-empty interior, we are going to define the tessellation with parameters $a, \Lambda$ as a random time process taking values in $\mathcal{M}$, with $a$ the time parameter. Let $B$ be Young's binary infinite tree, where the leaves are labeled so that the top leaf is attributed to number $1$, and the $k$-labeled leaf has the leaves $2k$ and $2k+1$ as daughter leaves. We now attach to each leave a pair $(\epsilon_{k}, \gamma_{k})$, where $(\epsilon_{k})_{k\geq 1}$ is a family of iid random exponential variables with parameter $\Lambda([W])$, and $(\gamma_{k})_{k\geq 1}$ are iid random hyperplanes with law $\frac{\Lambda([W]\cap .)}{\Lambda[W]}$.\\
 The tessellation is defined as a process of cell division, where each cell is identified with a leaf of the tree, and the first cell is the window $C_1=W$ itself. We describe it in term of a birth and death process, where a cell dies when it is divided in two daughter cells. The death times $(d_k)$, birth times $(b_k)$, and daughter-cells of $(C_k)_{\geq 1}$ are defined recursively in the following way ($\lfloor s \rfloor$ is the integer part of $s$):
$$\left\{
\begin{array}{rl}
b_1&=0\\
b_k&=d_{\lfloor k/2 \rfloor}\\
d_{k}&=b_k+\epsilon_{k}\\
C_{2k}&=C_k \cap \gamma_k^-\\
C_{2k+1}&=C_k \cap \gamma_k^+\\
\end{array}
\right.$$
We then define
\begin{eqnarray*}
R_{a,W}&=&\{C_k~|~b_k\leq a, C_k\neq \emptyset\},\\
Y_{a,W}&=&\bigcup_{b_k \leq a} \partial C_k.
\end{eqnarray*}
We call in the sequel $P_{a,W}$ the law of this tessellation. Note that $\Lambda$ is an implicit parameter of the model.\\
We also write in short
\begin{eqnarray*}
T_{Y_{a,W}}:=T_{a,W},~ U_{a,W}=1-T_{a,W}.
\end{eqnarray*}

\paragraph{Remarks}

\begin{enumerate}
\item The only non-intuitive feature of this construction is that $\gamma_{k}$ might not be in $[C_{k}]$, which would make one daughter cell be $\emptyset$ and the other be $C_{k}$ itself. There is no theoretical objection, but one can remedy this by, instead of choosing an iid family $(\epsilon_{k})$, attach independantly to each cell $C_{k}$ the death rate $\Lambda([C_{k}])$ and an hyperplane drawn from $\frac{\Lambda([C_{k}]\cap .)}{\Lambda[C_{k}]}$, so that the law is not modified and each hyperplane indeed hits the cell to which it is attached.
\item  When $\Lambda$ is isotropic, the death rate of a cell $C$,  $\Lambda([C])$, is proportionnal to the perimeter of $C$.
\item We mentioned in the introduction a construction from a Poisson hyperplane process. Let $\Theta$ be a Poisson process on $\mathbb{R}^+ \times \mathcal{H}$ with intensity $\Lambda([W])\lambda^+ \otimes \frac{\Lambda([W]\cap .)}{\Lambda[W]}$. Since a.s, for all $ t \geq 0,~  \mathop{\mathrm{card}}((\{t\}\times \mathcal{H}) \cap \Theta)\leq 1$, we can define a random sequence $(\tau_k, \gamma_k)_{k\geq 1}$ such that, almost surely, for all $ k \geq 1, \tau_k<\tau_{k+1}$, $\Theta=\{(\tau_k,\gamma_k)~|~k \geq 1\}$ and $(\tau_{k+1}-\tau_k, \gamma_k)$ satisfy the hpotheses of the previous construction. 
\end{enumerate}

\paragraph{Consistency}~\\

Nagel and Weiss \cite{NW1} have established the consistency property:
\begin{T}
If  $W \subset W'$ are two compact sets with non-empty interiors, then
\begin{eqnarray*}
Y_{a,W} \stackrel {(d)} {=}( Y_{a,W'} \cap W) \cup \partial W .
\end{eqnarray*}
\end{T}~\\
In particular, for any compact set $K \subset \mathop{\mathrm{int}}(W)$,
\begin{eqnarray*}
T_{a,W}(K)=T_{a,W'}(K),
\end{eqnarray*}
and we can define on $\mathcal{K}$
\begin{eqnarray*}
T_a(K):&=&T_{a,W}(K)\\
U_a(K):&=&U_{a,W}(K) \textrm{  ~~  for any $W$ containing $K$ in its interior.}
\end{eqnarray*}
Due to Choquet's Theorem, this allows us to define a tessellation $Y_a$ on the whole plane, with capacity functional $T_a$, such that
\begin{eqnarray*}
\A ~W \in \mathcal{K},~( Y_a \cap W) \cup \partial W \stackrel{(d)}{=} Y_{a,W},
\end{eqnarray*}
which is called the STIT tessellation with parameters $a$ and $\Lambda$. We call $P_a$ its law on $\mathcal{M}$ or, equivalently, on $\mathcal{M}_c$. Mecke, Nagel and Weiss \cite{MNW2} also provided an explicit global construction on $E$.

\paragraph{Iteration}~\\

Let us now define briefly the nesting operation. A more elegant definition, with the help of marked hyperplane Poisson processes, is given in \cite{MNW1}. Given two random tessellations $Y$ and $Y'$ having laws $Q$ and $Q'$, we define the law $Q  \boxplus Q'$ as that of the tessellation $Z$ obtained in the following way:\\
Let us take a copy of $Y$ and an independant iid family $(Y'_{q})_{q \geq 1}$ of copies of $Y'$. Let us put on the cells of $Y$ a measurable labelling independant of the $Y'_q$, $Y=\bigcup_{q \geq 1} \partial C_{q}$, and let us define $Z$ as the random closed set:
\begin{eqnarray*}
Z=Y \cup \left(\bigcup_{q \geq 1} C_{q} \cap Y'_{q}\right)
\end{eqnarray*}
so that $Z$ is obtained by nesting an independant copy of $Y'$ in each cell of $Y$. We note $Z= Y \circledcirc (Y_q')_{q \geq 1}$.\\
 A random tessellation is said to be stable under iteration  iff its law $P$ satisfies
 $$P\left(\frac{1}{2}.\right)=P \boxplus P.$$
 It has been shown in \cite{NW1} that STIT tessellations are the only random tessellations stable under iteration. Moreover, the class of STIT tessellations is attractive for the iteration in the following sense. Let us take a random stationary tessellation $Z$ with law $Q$, such that $a$ is the meanlength of $Z$ and $\Lambda$ is its directional distribution. We define the sequence of laws of tessellations $(Q_{n})_{n\geq 1}$ by:
 \begin{eqnarray*}
Q_1&=&Q,\\
Q_{n+1}&=&Q_{n}\boxplus Q_{1}.
\end{eqnarray*}
Hence, the law $Q_n$ is that of $Z$ iterated $n$ times with itself. It is shown in $\cite{NW1}$ that the normalized probability $Q_n(\frac{1}{n}.)$  converges weakly to the law of the STIT tessellation with parameters $a$ and $\Lambda$. So, the STIT tessellation is really the central object regarding this nesting operation. In a way it is not surprising since it is itself constructed in each compact window as some kind of iteration of the tessellation generated by a single hyperplane.\\

\section{Strong mixing property}

 In this part we will prove relation (\ref{EquMix}), and give a rate of decay close to the optimal for compact sets.\\
 Here $\Lambda=\lambda^+ \otimes \nu$ is a measure on $\mathcal{H}$ satysfyings assumption  (\ref{ConNu2}), and $a$ is a parameter strictly positive.\\
 In this section the probability space is $(\mathcal{M},\mathcal{B}(\mathcal{M}), P_a)$, and $Y_a$ is a random variable with law $P_a$.\\
Consider the set $\mathcal{T}=\{T_{h}~|~h \in \mathbb{R}^d\}$ of all translations, seen as operators on $\mathcal{F}$. The stationarity of a random closed set $X$ with law $P$ ensures that for all $T \in \mathcal{T}$, $P$ is invariant under $T$. Now we call $\mathcal{T}$-invariant set of $\mathcal{B}(\mathcal{F})$ every Borel set $\mathcal{C}$ of the $\sigma$-algebra $\mathcal{B}(\mathcal{F})$ such that $$\A ~T \in \mathcal{T}, ~T \mathcal{C}=\mathcal{C}.$$ For instance $\mathcal{K}, \mathcal{H}, \mathcal{M}$ or ``the class of all tessellations having a cube as one of their cells'' are invariant sets. Given a stationary law $P$ on $\mathcal{F}$, the dynamical system $(\mathcal{F},\mathcal{B}(\mathcal{F}),P,\mathcal{T})$ is said to be ergodic if every $\mathcal{T}$-invariant set has probability $0$ or $1$, and strongly mixing if for all $\mathcal{C}, \mathcal{C}'$ in $\mathcal{F}$, 
\begin{eqnarray}
\label{EquMixGal}
{P}(\mathcal{C},T_{h} \mathcal{C}' )\xrightarrow[ ||h||\to \infty]{h \in E} P(\mathcal{C})P(\mathcal{C}').
\end{eqnarray}
When $X$ is a random closed set with law $P$, we simply speak of the stationarity, the ergodicity or the mixing property for $X$. The mixing property ensures the asymptotic independance of relatively far away sets. Since any $\mathcal{T}$-invariant set $\mathcal{C}$ satisfies  ${P}(\mathcal{C},T_{h} \mathcal{C})=P(\mathcal{C})$, the strong mixing property implies ergodicity. \\
According to Lemma 9.3.1 in \cite{SW}, it suffices to show (\ref{EquMixGal}) for sets $\mathcal{C},\mathcal{C}'$ drawn from a semi-algebra $\mathcal{A}$ generating $\mathcal{C}$. That is why, in order to obtain the strong mixing property for $Y$, it is sufficient to show the following theorem:

 \begin{T}

 For all $A,B \in \mathcal{B},$

 \begin{eqnarray*}
 P_a(\mathcal{F}^{A} \cap \mathcal{F}^{T_h B})\xrightarrow[ ||h||\to \infty]{h \in E} P_a(\mathcal{F}^{A})P_a(\mathcal{F}^{B}).
 \end{eqnarray*}
 \end{T}

Note that ${T}_{h} \mathcal{F}^B=\mathcal{F}^{T^{-h}B}$, but we state the theorem with $h$ instead of $-h$ for more simplicity, since the role of $h$ is symmetric. Since the $(\mathcal{F}^{K})_{K \in \mathcal{K}}$ also generate the Borel $\sigma$-algebra of the Fell topology, still using Lemma 9.3.1 in \cite{SW}, it suffices to show it for $A,B \in \mathcal{K}$, and the proof is given in section \ref{SecPrf}. We are able to indicate the optimal rate of decay, for that we first need to introduce a function.

\begin{Ps}
\label{ProFunDir}

Given $u \in \mathcal{S}^{d-1}$, the measure of all hyperplanes separating $0$ and $u$ has the following representation:

\begin{eqnarray*}
\zeta(u):=\Lambda([0|u])=\Lambda([0,u])=\frac{1}{2} \int_{\mathcal{S}^{d-1}}\nu(\d v) |\langle u,v\rangle|,
\end{eqnarray*}
and so
\begin{eqnarray*}
\A ~u, v \in \mathcal{S}^{d-1},~ |\zeta(u)-\zeta(v)| \leq \frac{1}{2}||u-v||.
\end{eqnarray*}
Furthermore, there exists $\kappa>0$ such that
\begin{eqnarray*}
\A ~u \in \mathcal{S}^{d-1}, \zeta(u)\geq \kappa.
\end{eqnarray*}

\end{Ps}
This result, shown in section \ref{SecPrf} along with others, is meaningful regarding the asymptotic behaviour of expression (\ref{EquMix}), given in the following theorem:

\begin{T}
\label{ThmRatMix}
 For all $ A,B \in \mathcal{K},$ There exists  a constant $ \chi(A,B,a)>0$, invariant to rigid motions of $A$ and $B$, such that,
\begin{eqnarray}
\label{EquRatCom}
\left|P_a(\mathcal{F}^{A} \cap \mathcal{F}^{T_h B})-P_a(\mathcal{F}^{A})P_a(\mathcal{F}^{B})\right|\leq \frac{\chi(A,B,a)}{||h||\zeta(h/||h||)}(1+o(1)).
\end{eqnarray}
Furthermore, when $A$ and $B$ each have a single connected component,
\begin{equation}
\label{EquRatCon}
\left|\frac{P_a(\mathcal{F}^A,\mathcal{F}^{T_h B})}{P_a(\mathcal{F}^A)P_a(\mathcal{F}^B)}-1\right|=\frac{1}{||h||\zeta(h/||h||)}(1+o(1)).
\end{equation}
\end{T}

The asymptotic dependence of $\mathcal{F}^A$ and $\mathcal{F}^{T_h B}$ is closely connected to the probability that the cell $W=\mathop{\textrm{conv}}(A \cup T_h B)$ is hit before time $a$ by a hyperplane that separates $A$ and $T_h B$, which dependency is asymptotically in $||h||\zeta(h/||h||)$, using the stationarity of $\Lambda$ (see lemma \ref{LmSepPts}).\\
According to (\ref{EquRatCon}), (\ref{EquRatCom}) is optimal in the sense that the inequality is an equality when $A$ and $B$ are connected sets. When they are not, we think that we are close to optimality, but the estimation is more difficult  because the capacity functional has a much more complicated expression, even with only two components (see \cite{NW1}). To deal with it, we have to estimate the lipschitziannity of the capacity functional as a temporal function of $a$. We obtain the following expression, worth to be noted:\\
Let $t>0$ be a real positive number, and $K$ a compact set. Then
\begin{equation}
\label{EquMajCap}
\textrm{There exists } \lambda_{K,a}>0,~~~
0 \leq T_{a+t}(K) -T_{a}(K)  \leq t \lambda_{K,a}.
\end{equation}

Furthermore, $\lambda_{K,a}$ is invariant to rigid motions of $K$.\\
If one had a non-zero linear lower bound in (\ref{EquMajCap}), one could probably, by a slight modification of the proof of theorem \ref{ThmRatMix}, obtain a lower bound of the left-hand term in (\ref{EquRatCom}) by one of the form of the right-hand term.

\section{Proofs}
\label{SecPrf}
\paragraph{Proof of (\ref{EquMajCap})}
\begin{Lm}
Let $0<t,a$ be real positive numbers, $K$ a compact set, and \\$\lambda_{K,a}=\Lambda([\mathop{\mathrm{conv}}(K)])(1+a\Lambda([\mathop{\mathrm{conv}}(K)]))(1-T_{a(K)})$. Then
\begin{eqnarray*}
0 \leq T_{a+t}(K) -T_{a}(K)  \leq t \lambda_{K,a} .
\end{eqnarray*}
\end{Lm}~\\
We use the notations of section \ref{SecSTIT}, with $W=\mathop{\mathrm{conv}}(K)$.
Let us call $N_a$ the number of times $\mathop{\mathop{\textrm{conv}}}(K)$ has been hit by a hyperplane up to time $a$:
\begin{eqnarray*}
N_a=\mathop{\mathrm{card}}\{~k~ |~\gamma_k \in [C_k]\}.
\end{eqnarray*}
 For $n \in \mathbb{N}$ we set $\pi_n(a)=\mathbb{P}(N_a=n)>0$. Almost surely, conditionaly to $(N_a=n)$, $R_{a,W}$ contains $n+1$ cells $C_1,...,C_{n+1}$ with non-empty interior. We call  $K_i:=K \cap C_i$.
We have
\begin{multline*}
T_{a+t}(K) -T_{a}(K)  =U_a(K)-U_{a+t}(K)=U_a(K) \mathbb{P}(K \cap Y_{a+t}\neq \emptyset|K \cap Y_a=\emptyset)\\
=U_a(K)\sum_{n \geq 0} \pi_n(a) \mathbb{P}(\Ex i, 1 \leq i \leq n+1, K_i \cap Y_{a+t}\neq\emptyset | N_a=n, K \cap Y_a=\emptyset )\\
\leq U_a(K)\sum_{n \geq 0} \pi_n(a) \sum_{i=1}^{n+1}\mathbb{P}(K_i \cap Y_{a+t}\neq\emptyset | N_a=n, K \cap Y_a=\emptyset )\\
= U_a(K)\sum_{n \geq 0} \pi_n(a) \sum_{i=1}^{n+1}\mathbb{P}(\mathop{\textrm{conv}}(K_i) \cap Y_{a+t}\neq\emptyset |N_a=n,  \mathop{\textrm{conv}}(K_i) \cap Y_a = \emptyset ).\\
\end{multline*}
At this point, for $\mathop{\mathop{\textrm{conv}}}(K_i)$, knowing that $K_i$ is in the interior of a cell at time $a$ amounts to know that $\mathop{\textrm{conv}}(K_i)$ has not been hit up to time $a$, and so, due to the lack of memory and the consistency property, 
\begin{eqnarray*}
\mathbb{P}(\mathop{\mathop{\textrm{conv}}}(K_i) \cap Y_{a+t}\neq\emptyset~ |~N_a=n,  \mathop{\textrm{conv}}(K_i) \cap Y_a = \emptyset) &=& 1-U_t(\mathop{\mathop{\textrm{conv}}}(K_i))\\&\leq & 1-U_t(\mathop{\textrm{conv}}(K))\\ 
&=&e^{-t \Lambda([\mathop{\textrm{conv}}(K)])}.
\end{eqnarray*}
This yields
\begin{eqnarray*}
U_{a}(K)-U_{a+t}(K)&\leq&U_a(K)\sum_{n \geq 0}\pi_n(a) \sum_{i=1}^{n+1} \left(1-e^{-t\Lambda([\mathop{\mathop{\textrm{conv}}}(K)])}\right)\\
&\leq&U_a(K)\sum_{n \geq 0}\pi_n(a) t\sum_{i=1}^{n+1}\Lambda([\mathop{\mathop{\textrm{conv}}}(K)])\\
&\leq&U_a(K)\sum_{n\geq 0}\pi_n(a)(n+1) t\Lambda([\mathop{\mathop{\textrm{conv}}}(K)])\\
&\leq&U_a(K) t\Lambda([\mathop{\textrm{conv}}(K)])(1+\mathbb{E}(N_a)).
\end{eqnarray*}
The number of hyperplanes involved in the cell division process with initial cell $\mathop{\mathop{\textrm{conv}}}(K)$ up to time $a$ is less or equal to a Poisson variable with parameter $a\Lambda([\mathop{\textrm{conv}}(K)])$ (remember that some hyperplanes fall outside the cell they were supposed to divide). In consequence we have
\begin{eqnarray*}
\mathbb{E}(N_a) \leq a\Lambda([\mathop{\mathop{\textrm{conv}}}(K)]),
\end{eqnarray*}
which concludes the proof.
\begin{flushright}
$\Box$
\end{flushright}

\paragraph{Proof of  Proposition \ref{ProFunDir}}~\\

For $u \in \Sd$ we have by stationarity
\begin{equation*}
\zeta(u)=\Lambda(]0,u[)=\Lambda([0,u])=\frac{1}{2}\Lambda([-u,u]) \text{   by stationarity}.
\end{equation*}
Now, given a hyperplane $\gamma$ at distance $r$ from the origin with exterior normal vector $v$, 
\begin{equation*}
\gamma \in [-u,u] \Leftrightarrow |\langle u,v \rangle| \leq r,
\end{equation*}
and so
\begin{eqnarray*}
\zeta(u)=\frac{1}{2} \int_{\Sd} |\langle u,v \rangle|\nu(\d v).\\
\text{For  } u\neq u' \in \Sd, |\zeta(u)-\zeta(u')| &\leq& \frac{1}{2}\int_{\Sd}  \left| ~|\langle u,v \rangle|- |\langle u',v \rangle|~ \right| \nu(\d v)\\
&\leq &  \frac{1}{2}\int_{\Sd}  |\langle u,v \rangle- \langle u',v \rangle| \nu( \d v)\\
& \leq & \frac{||u-u'||}{2}\int_{\Sd} \nu(\d v)\\
& \leq &  \frac{||u-u'||}{2}\nu(\Sd).
\end{eqnarray*}
Let $u \in \Sd$. Thanks to assumption (\ref{ConNu2}), $\nu$ is not concentrated on $u^{\perp}$, hence $ \int_{\Sd} |\langle u,v \rangle|\nu(\d u) >0$. Since $\Sd$ is compact, we have the existence of $\kappa>0$ such that $\zeta(u)>\kappa$.

\begin{flushright}
$\Box$
\end{flushright}

\paragraph{Proof of Theorem \ref{ThmRatMix}}~\\

We are first going to establish some inequalities for $A,B$ any compact sets, and then add a drift $h$ to give an upper bound when expressions become too complicated.

Let $W=\mathop{\mathop{\textrm{conv}}}(A \cup B)$.\\
The key is to consider the tessellation valued time process defined in section \ref{SecSTIT},  $(Y_{t,W})_{t\geq 0}$. We keep the notations $\epsilon_{1}=\inf\{t~|~Y_{0,W}\neq Y_{t,W}\}$ and $\gamma_{1}$ is the first hyperplane dividing  $C_1= W$.\\
We introduce the event

\begin{eqnarray*}
\Gamma_{A,B}=(\gamma_{1}\in [A|B], \epsilon_1\leq a).
\end{eqnarray*} 
When $A$ and $B$ are far away from each other, $\Gamma_{A,B}$ is very likely to happen, and after $\Gamma_{A,B}$ occurs the behaviours of $A$ and $B$ are independants because they are in disjoint cells.\\
We first have to estimate
\begin{eqnarray*}
| P_a(\mathcal{F}^A,\mathcal{F}^B,\Gamma_{A,B})-P_a(\mathcal{F}^A,\mathcal{F}^B)|& \leq & P_a(\Gamma_{A,B}^c)\\
&\leq & e^{-a\Lambda([W])}.\\
\end{eqnarray*}
In consequence, we can show formulae (\ref{EquRatCom}) and (\ref{EquRatCon}) with $P_a(\mathcal{F}^A,\mathcal{F}^{T_h B},\Gamma_{A,T_h B})$ instead of $P_a(\mathcal{F}^A,\mathcal{F}^{T_h B})$, because their difference is negligible with respect to the expected decay rate.
Rigorously this gives:
\begin{multline*}
P_a(\mathcal{F}^A,\mathcal{F}^B,\Gamma_{A,B})=\int_{0}^a P_a(\epsilon_{1} \in \d t, \gamma_{1} \in [A|B], \mathcal{F}^A, \mathcal{F}^B)\\
=\int_{0}^a P_a(\epsilon_{1} \in \d t, \gamma_{1}  \in [A|B])P_a(\mathcal{F}^A, \mathcal{F}^B | \epsilon_{1}=t, \gamma_{1}  \in [A|B]) .
\end{multline*}
When $\epsilon_{1}=t$ and $\gamma_{1} \in [A|B]$, it means that $A$ and $B$ have not been hit up to time $t$ and are both still contained in a cell (but the cell encapsulating $A$ is different from that of $B$). Thanks to the consistency property, and the independent behaviour of distinct cells after their birth, we have

\begin{eqnarray*}
P_a(\mathcal{F}^A, \mathcal{F}^B | \epsilon_{1}=t, \gamma_{1}  \in [A|B])=U_{a-t}(A)U_{a-t}(B).
\end{eqnarray*}
Since the sequence $(\gamma_{k})$ is independent of $(\epsilon_{k})$, we have 
\begin{eqnarray*}
P_a(\epsilon_{1} \in \d t, \gamma_{1}  \in [A|B])&=& P_a(\epsilon_{1} \in \d t) P_a(\gamma_{1}  \in [A|B])\\
&=&\Lambda([W])e^{-t \Lambda([W])}\d t \frac{\Lambda([A|B])}{\Lambda([W])},
\end{eqnarray*}
and one finally has
\begin{eqnarray}
\label{EquExaMix}
P_a(\mathcal{F}^A,\mathcal{F}^B,\Gamma_{A,B})&=&\Lambda([A|B])\int_{0}^a e^{-t \Lambda([W])} U_{a-t}(A) U_{a-t} (B) \d t,
\end{eqnarray}
so that
\begin{multline*}
|P_a(\mathcal{F}^A,\mathcal{F}^B,\Gamma_{A,B})-P_a(\mathcal{F}^A)P_a(\mathcal{F}^B)|\\
 = \left| \Lambda([A|B])\int_{0}^a e^{-t \Lambda([W])} U_{a-t}(A) U_{a-t} (B) \d t-\right.\\
\left. U_a(A)U_a(B)\int_{0}^{\infty} e^{-t\Lambda([W])} \Lambda([W])\d t \right|
\end{multline*}
\begin{multline}
\label{EquMajMix}
\shoveleft \leq  {\Lambda([A|B])}\int_{0}^{\infty} e^{-t \Lambda([W])}\left| U_{a-t}(A) U_{a-t} (B)-
U_a(A)U_a(B)\right|  \d t+\\
\shoveright{ +\left|\frac{\Lambda([A|B])}{\Lambda([W])}-1\right|\int_{0}^{\infty} e^{-t \Lambda([W])} U_{a}(A) U_{a} (B) \Lambda([W])\d t+}. \\
+\Lambda([A|B]) \int_a^{\infty} e^{-t\Lambda([W])} \d t.
\end{multline}
Now it is time to consider the drift of $B$, in order to make a proper upper bound estimate. We will use  $T_h B$ instead of $B$ in the previous formulae and use the following facts.

\begin{Lm}Let $A, B \in \mathcal{K}$ and $ h$ be an element of $\mathbb{R}^d$. We set \\$W_h=\mathop{\mathrm{conv}}(A,T_h B)$. Then
\begin{eqnarray}
\int_{0}^{||h||} e^{-u} \d u&=&o(\frac{1}{||h||}),\\
\Lambda([W_{h}])&=& \Lambda([A|T_h B])(1+o(1)),\\
\Lambda([W_{h}])&=&\varphi(h/||h||)||h||(1+o(1)).
\end{eqnarray}
\end{Lm}
\begin{proof}
Let us prove the two latter formulas:
Let $\alpha \in A, ~\beta \in B$.
 If $\gamma$ is a hyperplane separating $\alpha$ and $\beta-h$ and hitting neither $\mathop{\textrm{conv}}(A)$ nor $\mathop{\textrm{conv}}(T_h B)$, it also separates $A$ and $T_h B$. Hence
\begin{multline*}
|\Lambda([A|T_h B])-\Lambda([\alpha,\beta-h])|\leq \Lambda([\mathop{\mathop{\textrm{conv}}}(A)])+\Lambda([\mathop{\mathop{\textrm{conv}}}(T_h B)]) \\
\leq \Lambda([\mathop{\mathop{\textrm{conv}}}(A)])+\Lambda([\mathop{\mathop{\textrm{conv}}}(B)] =o(||h||).
\end{multline*}
For similar reasons,
\begin{eqnarray*}
0 \leq \Lambda(W_{h})-\Lambda([A|T_h B) \leq \Lambda([\mathop{\mathop{\textrm{conv}}}(A)])+\Lambda([\mathop{\mathop{\textrm{conv}}}(B)]) = o(||h||),
\end{eqnarray*}
and we only need to prove

\begin{eqnarray*}
\Lambda([\alpha, \beta-h])&=&\varphi(h/||h||)||h||(1+o(1)),
\end{eqnarray*}
which is given by the following result.
\begin{Lm}
\label{LmSepPts}
If $\beta-\alpha=||\beta-\alpha||u$, with $u \in \mathcal{S}^{d-1}$,
$$\Lambda([\alpha|\beta])=||\beta-\alpha||\zeta(u).$$
\end{Lm}
\begin{proof}

It suffices to show that for all $\epsilon>0, u \in \mathcal{S}^{d-1}, n \in \mathbb{N}$,

\begin{eqnarray*}
  \Lambda([0 |(n+1)\epsilon u]) = \Lambda([0|n\epsilon u]) + \Lambda([0|\epsilon u]) ,
\end{eqnarray*}
because, since $\Lambda$ is stationnary and locally finite,  $\A ~x \in \mathbb{R}^d, \Lambda([x])=0$ and then we will be able to obtain the result by induction and approximation.\\
Since $ [0|(n+1)\epsilon u]$ is the disjoint union of $[0|n \epsilon u], [n \epsilon u|(n+1)\epsilon u] $ and $[n\epsilon u]$, we have

\begin{eqnarray*}
  \Lambda([0|(n+1)\epsilon u]) = \Lambda([0|n\epsilon u]) + \Lambda([n \epsilon u|(n+1)\epsilon u]) 
\end{eqnarray*}
and the result follows by stationnarity of $\Lambda$.\\

\end{proof}
Finally we can end the proof of the lemma:
\begin{eqnarray*}
\Lambda([W_{h}])=||\alpha-\beta+h||\zeta(\frac{\alpha-\beta+h}{||\alpha-\beta+h||})(1+o(1))=||h||\zeta(h/||h||)(1+o(1))
\end{eqnarray*}
since $\zeta$ is continuous.

\end{proof}

Let us define
\begin{eqnarray*}
\delta(A,B,a)&=&\lambda_{A,a}+\lambda_{B,a}.\\
&\geq  & \lambda_{A,a-t}+\lambda_{B,a-t}~ \textrm{ for all $t>0$.}
\end{eqnarray*}
We set $$\chi(A,B,a)=\delta(A,B,a)+\Lambda([\mathop{\textrm{conv}}(A)])+\Lambda([\mathop{\textrm{conv}}(B)]).$$
Using (\ref{EquMajCap}) and (\ref{EquMajMix}), and the previous lemma, we get
\begin{multline*}
|P_a(\mathcal{F}^A,\mathcal{F}^{T_h B},\Gamma_{A,T_h B})-P_a(\mathcal{F}^A)P_a(\mathcal{F}^B)|\\
\shoveright{ \leq  \Lambda([A|T_h B]) \delta(A,B,a) \int_{0}^{\infty}e^{-t\Lambda([W_{h}])} t \d t + (1-\frac{\Lambda([A|T_h B])}{\Lambda([W_{h}])})+o(\frac{1}{||h||})}\\
\shoveright{\leq \frac{\delta(A,B,a) \Lambda([A|T_h B])}{\Lambda([W_{h}])^2}+\frac{\Lambda([\mathop{\textrm{conv}}(A)])+\Lambda([\mathop{\textrm{conv}}(T_h B)])}{\Lambda([W_{h}])}+o(\frac{1}{||h||})}\\
\leq  \frac{\chi(A,B,a)}{\varphi(h/||h||)||h||}(1+o(1)),
\end{multline*}
and formula (\ref{EquRatCom}) is proved.\\

When $A$ and $B$ are connected, thanks to the simple expression (see \cite{NW1})
 $$U_u(A)=e^{-u\Lambda([A])},$$
 using (\ref{EquExaMix}) we have the exact expression
\begin{multline*}
|P_a(\mathcal{F}^A,\mathcal{F}^{T_h B},\Gamma_{A,T_h B})-P_a(\mathcal{F}^A)P_a(\mathcal{F}^B)| \\
=P_a(\mathcal{F}^A)P_a(\mathcal{F}^B) \int_0^a e^{-t(\Lambda(W_{h})-\Lambda(A)-\Lambda(B))} \d t\\
=  P_a(\mathcal{F}^A)P_a(\mathcal{F}^B) \left(\frac{1}{\Lambda(W_{h})-\Lambda(A)-\Lambda(B)}\right)\int_{0}^{\infty}e^{-u} \d u (1+o(1)).\\
\end{multline*}
Hence
\begin{equation}
\left|\frac{P_a(\mathcal{F}^A,\mathcal{F}^{T_h B},\Gamma_{A,T_h B})}{P_a(\mathcal{F}^A)P_a(\mathcal{F}^{B})}-1\right|=\frac{1}{||h||\varphi(h/||h||)}(1+o(1)),
\end{equation}

and (\ref{EquRatCon}) is proved.

\begin{flushright}
$\Box$
\end{flushright}

\section{Acknowledgements} 
I am grateful to Werner Nagel, for useful discussion and comments. I also thank Youri Davydov, my PhD advisor, and Christoph Thaele, for their careful reading and remarks,  and the workgroup of stochastic geometry from Lab. Paul Painlev\'e, Univ. Lille 1.

\end{document}